\documentclass[11pt]{amsart}
\usepackage{amsfonts}
\usepackage{amsmath}
\usepackage{amssymb}
\usepackage{graphicx}
\usepackage{amsthm,amscd}
\usepackage{calc}
\usepackage{amsthm}

\newtheorem{theorem}{Theorem}[section]

\newtheorem{corollary}{Corollary}[section]

\newtheorem{proposition}{Proposition}[section]

\numberwithin{equation}{section}
\theoremstyle{remark}

\renewcommand{\bar}{\overline}
\newcommand{\eps}{\varepsilon}
\newcommand{\pa}{\partial}
\renewcommand{\phi}{\varphi}
\newcommand{\wt}{\widetilde}

\newcommand{\ka}{K\"ahler }
\newcommand{\kr}{K\"ahler-Ricci }

\newcommand{\C}{{\mathbb C}}

\newcommand{\R}{{\mathbb R}}

\newcommand{\M}{{\mathcal M}}

\newcommand{\XX}{{\mathfrak X}}

\newcommand{\ke}{K\"ahler-Einstein }

\newcommand{\h}{{\mathbb H}}

\newcommand{\g}{{\mathfrak g}}

\newcommand{\lb}{\left (}
\newcommand{\rb}{\right )}
\newcommand{\lsb}{\left [}
\newcommand{\rsb}{\right ]}
\newcommand{\lfb}{\left \{}
\newcommand{\rfb}{\right \}}

\newcommand{\dd}{\text{div}}
\newcommand{\ga}{\alpha}
\newcommand{\gb}{\beta}

\newcommand{\rc}{\text{Ric}}
\newcommand{\ji}{\mathcal J^{int}}
\newcommand{\wpm}{Weil-Petersson }
\newcommand{\ii}{\sqrt{-1}}
\newcommand{\au}{\text{Aut}_0}
\newcommand{\tr}{\text{Tr}}

\title[]{A Weil-Petersson Type Metric on the Space of Fano K\"ahler-Ricci Solitons}

\author{Huai-Dong Cao$^1$, Xiaofeng Sun$^1$ and Yingying Zhang$^2$}
\address{Department of Mathematics,  Lehigh University,
Bethlehem, PA 18015, USA}
\email{huc2@lehigh.edu}
\address{Department of Mathematics,  Lehigh University,
Bethlehem, PA 18015, USA}
\email{xis205@lehigh.edu}
\address{Yau mathematical Sciences Center, Tsinghua University, Beijing, 100804, China}
\email{yingyzhang@tsinghua.edu.cn}

\dedicatory{Dedicated to Professor Peter Li  on  the occasion of his 70th birthday}

\thanks{$^1$Research partially supported by Simons Foundation Collaboration Grants.}
\thanks{$^2$Research supported in part by NSFC Grant No. 12141101}

\begin{document}

\begin{abstract} In this paper we define a Weil-Petersson type metric on the space of shrinking K\"ahler-Ricci solitons and prove a necessary and sufficient condition on when it is independent of the choices of K\"ahler-Ricci soliton metrics. We also show that the Weil-Petersson metric is K\"ahler when it defines a metric on the Kuranishi space of small deformations of K\"ahler-Ricci solitons.  Finally, we establish the first and second order deformation of Fano K\"ahler-Ricci solitons and show that, essentially, the first effective term in deforming K\"ahler-Ricci solitons leads to the Weil-Petersson metric.
\end{abstract}

\maketitle

\section{Introduction}\label{intr}

The \wpm metric is a powerful tool in understanding the geometry of moduli spaces of various geometric objects equipped with certain canonical metrics. It was first introduced by Weil \cite{weil} in the late 1950s by using the Petersson pairing on modular forms and it is a natural $L^2$ metric on the moduli space $\M_g$ of hyperbolic Riemann surfaces of genus $g\ge 2$. More precisely, for any point $p$ in $\M_g$ represented by a Riemann surface $\Sigma$,  %the Kodaira-Spencer theory implies
it is well-known that the holomorphic tangent space of $\M_g$ at $p$ can be identified with $H^{0,1}\lb \Sigma, T^{1,0}\Sigma\rb$, and that the dual space $\Omega_p^{1,0}\M_g$ can be identifies with $H^0\lb K_\Sigma^2\rb$, the space of holomorphic quadratic differentials on $\Sigma$. The Weil-Petersson co-metric, as defined by Weil, is given by $\Vert \eta\Vert_{_{WP}}^2=\int_\Sigma |\eta|_\lambda^2 dV_\lambda,$
where $\lambda$ is a hyperbolic metric on $\Sigma$ and $\eta\in H^0\lb K_\Sigma^2\rb$. It follows that the \wpm metric on $H^{0,1}\lb \Sigma, T^{1,0}\Sigma\rb$ is defined by
\[
\Vert \phi\Vert_{_{WP}}^2=\int_\Sigma |\phi|_\lambda^2 dV_\lambda
\]
for any harmonic Beltrami differential $\phi$ with respect to $\lambda$.

The \wpm metric on $\M_g$ is well-understood.  In the 1960s, Ahlfors \cite{ahlforsa, ahlforsb} showed that the \wpm metric on $\M_g$ is \ka and has negative curvature; see also Wolpert's work \cite{wolpert}. Chu \cite{chu} and Masur \cite{mas} computed the asymptotics of the \wpm metric and showed that the \wpm metric is incomplete and has finite volume.

Subsequently, the \wpm metric has been generalized to study the geometry of moduli spaces of higher-dimensional varieties and vector bundles which admit certain canonical metrics.
Koiso \cite{koiso83} first studied the Weil-Petersson metric on the moduli space of \ke manifolds and proved that it is K\"ahler. Royden \cite {ro1974} and Siu \cite{siu} developed the canonical lifting method to compute the curvature of the \wpm metric on the moduli spaces of \ke manifolds of general type. Their method was adopted by Nannicini \cite{nan} to study  the \wpm metric on the moduli spaces of polarized Calabi-Yau manifolds. Candelas, Green and H\"ubsch \cite{cgb88} showed that the \wpm metric on the moduli spaces of polarized Calabi-Yau manifolds is the curvature of the $L^2$ metric on the first Hodge bundle and its curvature can be derived from the work of Todorov \cite{tod}.  In these cases, the \wpm metric on the parameter space of a family of \ke manifolds is defined as the $L^2$ metric on the space of harmonic Beltrami differentials on each fiber of the family with respect to fiberwise \ke metrics. The lifting method was later refined by Schumacher \cite{sch} and, together with Fujiki, they introduced the generalized \wpm metric on the (coarse) moduli space of csck metrics  \cite{fs1990}.

The \wpm metric is well-defined as long as the \ke metric on each $X_t$ is unique. This is the case when the fibers are of general type, Calabi-Yau type, or Fano type with discrete automorphism groups. However, when the fibers are Fano \ke manifolds with non-discrete automorphism groups, there is a family of \ke metrics on each fiber.  In this situation, the $L^2$ inner product and harmonic representatives of the Kodaira-Spencer classes may vary according to the choices of \ke metrics. Thus, it is not a priori clear that the \wpm metric is well-defined.  In fact, a necessary and sufficient condition was provided  in \cite{csyz2021} to guarantee the well-definedness of the Weil-Petersson metric in the Fano \ke case (see also Theorem~\ref{kewell}). There is also a canonical \wpm metric on the (complexified) \ka moduli spaces which was defined in a similar fashion \cite{wilson04, trwil11}. Such Weil-Petersson metrics play an important role in mirror symmetry.

In this paper, we study the Fano \kr soliton case and extend the work in \cite{csyz2021} to define and investigate a \wpm type metric on the Kuranishi space of small deformations of Fano \kr solitons. This metric is the natural $L^2$ metric on the space of $f$-twisted harmonic Beltrami differentials, where $f$ is the Ricci potential of a \kr soliton metric.  We also provide a necessary and sufficient condition for the \wpm metric  to be well-defined  in the Fano \kr soliton case; see Theorem~\ref{wd}.  It turns out, as shown in Theorem~\ref{2nd}, that this \wpm metric is closely tied to the variation of \kr solitons as the underlying complex structures vary.

We remark that Donaldson, in \cite{do2015}, defined a new \ka metric on the space $\mathcal J$ of complex structures that are compatible with a given symplectic structure $\omega$ on the background smooth manifold of a Fano manifold. In particular, the moment map of the action of the symplectomorphism group on $\mathcal J$, with respect to this new \ka structure, is given by the volume form of $\omega$ twisted by the Ricci potential. It follows that the zeros of the moment map correspond to \ke metrics, while Fano \kr solitons are just the critical points of the $H$-functional defined in \cite{he}, which can be view as the $``$norm square$"$ of the moment map. In this situation, the \wpm metric we defined can be viewed as the restriction of Donaldson's metric on the Kuranishi space of Fano \kr solitons when the twisted Kuranishi-divergence gauge is imposed.

The paper is organized as follows. In Section \ref{kdga}, we review the twisted Kuranishi-divergence gauge for general Fano manifolds that was established by the third author in her 2014 thesis \cite{zhang}. In Section \ref{wpmc}, we first recall small deformation of Fano \ke manifolds and the \wpm metric in the Fano \ke setting,  and then explain Donaldson's construction of his new metric on $\mathcal J$ and the associated moment map. We then define the \wpm metric on the space of Fano \kr solitons. To get a useful metric on the moduli space of Fano \kr solitons, the \wpm metric we introduced has to be independent of the choices of \kr soliton metrics. We will prove a  sufficient and necessary condition on such independence and show that the \wpm metric is \ka when it defines a metric on the Kuranishi space. Finally, we establish the first and second order deformation of Fano \kr solitons and show that essentially the first effective term in deforming \kr solitons leads to the \wpm metric.

\section{The Twisted Kuranishi-Divergence Gauge}\label{kdga}

Let $\lb M^n, \omega_g\rb$ be a Fano manifold of dimension $n$ such that $\lsb \omega_g\rsb=2\pi c_1\lb M\rb$ with Ricci potential $f\in C^\infty(M, \R)$:
\[
\rc\lb\omega_g\rb=\omega_g+\ii\pa\bar\pa f.
\]
We normalize $f$ so that $\int_M e^f dV_g=\int_M \lb 2\pi c_1\lb M\rb\rb^n$, where \[dV_g=\omega_g^{[n]}:=\frac {\omega_g^n}{n!}\]
is the volume form of $\omega_g$.
Note that, by definition, $\omega_g$ is a (shrinking) \kr soliton if the vector field $X=\nabla^{1,0} f\in A^0\lb M, T^{1,0}M\rb$ is holomorphic, namely if $X\in H^0\lb M, T^{1,0}M\rb$. We refer the reader to \cite{cao, koiso90, wangzhu} for examples of Fano \kr solitons, and to \cite{cao2010} for a basic overview of Ricci solitons.

The twisted Laplacian on smooth functions is defined by
\[
\Delta_f u=\Delta u+\bar X\lb u\rb=g^{i\bar j}\frac{\pa^2 u}{\pa z_i\pa\bar{z_j}}+g^{i\bar j}\frac{\pa f}{\pa z_i}\frac{\pa u}{\pa\bar{z_j}}, \qquad u\in C^\infty \lb M, \C\rb
\]
% for any $u\in C^\infty \lb M, \C\rb$
where $\lb z_1, \dots, z_n \rb$ is any system of local holomorphic coordinates on $M$. We also denote by  $\bar\pa^{*}_f$ the twisted adjoint operator of $\bar\pa$ with respect to the weighted volume form $e^fdV_g$. It is natural in the current situation since
\[
\rc\lb e^fdV_g\rb=\omega_g.
\]
The twisted Laplacian $\Delta_f$ is not a real operator and we define its conjugate operator $\bar\Delta_f$ by
\[
\bar\Delta_f u=\bar{\Delta_f\bar u}=\Delta u+X\lb u\rb.
\]
Similarly, we can define the twisted divergence operator $\dd_f$ on vector fields by
 \[\dd_f\lb V\rb=\dd (V)+V(f).\]

Although $\Delta_f$ is not a real operator, it is self-adjoint with respect to the Hermitian inner product
\[
(u,v)=\int_M u\bar v\ e^fdV_g
\]
and thus its eigenvalues are real. Furthermore, the eigenvalues are all positive (except the zero eigenvalue when acting on constant functions). Let $\lambda_1$ be the smallest positive eigenvalue of $\Delta_f$ and denote by $\Lambda_f^1$, in case $\lambda=1$ is an eigenvalue,  the eigenspace of eigenvalue $1$:
\[
\Lambda_f^1 =\lfb u\in C^\infty\lb M,\C\rb\mid \lb 1+\Delta_f\rb u=0\rfb.
\]

Extending the Lichnerowicz theorem, Futaki \cite{futaki} showed that

\begin{theorem}\label{e1} Let $\lb M, \omega_g\rb$ be a Fano manifold with Ricci potential $f$. Then
$\lambda_1\equiv \lambda_1 (\Delta_f)\geq 1$. Moreover, if $H^0\lb M, T^{1,0}M\rb\ne 0$ then $\lambda_1=1$ and $H^0\lb M, T^{1,0}M\rb\cong \Lambda_f^1.$
In particular,
\[
\dd_f: H^0\lb M, T^{1,0}M\rb\to \Lambda_f^1 \ \ \ \text{and}\ \ \ \nabla^{1,0}:\Lambda_f^1\to H^0\lb M, T^{1,0}M\rb
\]
 are linear isomorphisms with $\nabla^{1,0}\circ \dd_f=-Id$.
\end{theorem}

Next we consider the Kodaira-Spencer-Kuranishi deformation theory of Fano manifolds. On the space of Beltrami differentials $A^{0,1}\lb M, T^{1,0}M\rb$ the twisted adjoint of the $\bar\pa$ operator is given by
\[
\bar\pa_f^* \ \! \phi=\bar\pa^*\phi-\bar X\lrcorner \ \! \phi
\]
and the twisted divergence operator is
\[
\dd_f \ \! \phi=\dd \ \!\phi+\phi\lrcorner \ \! \pa f
\]
for any $\phi\in A^{0,1}\lb M, T^{1,0}M\rb$. In the following, we will use $\h_f$ to denote the spaces of $f$-harmonic forms or the $L^2$ projection, with respect to the metric $g$ and the weighted volume form $e^fdV_g$, to such harmonic spaces. For example
\[
\h_f^{0,1}\lb M, T^{1,0}M\rb=\lfb \phi\in A^{0,1}\lb M, T^{1,0}M\rb\big | \ \bar\pa \phi=0,\ \bar\pa_f^* \phi=0\rfb.
\]
Note that, by the Kodaira vanishing theorem and the Serre duality, we have
\[
H^{0,k}\lb M, T^{1,0}M\rb=0 \ \ \ \text{for all}\ \  k\geq 2,
\]
thus the deformation of the complex structures on $M$ is unobstructed.

In 2014, the third author in her thesis \cite{zhang} obtained the following result on the deformation of complex structures with respect to the twisted Kuranishi gauge.

\begin{theorem}  \label{twku}
Let $\lb M_0,\omega_0\rb$ be a Fano manifold, with $\omega_0$ in the class $2\pi c_1\lb M_0\rb$, and let $f_0$ be the normalized Ricci potential. Let $m=h^{0,1}\lb M_0, T^{1,0}M_0\rb$ and let $0\in B\subset \C^m$ be the open ball of radius $\eps$ with coordinates $t=\lb t_1,\cdots,t_m\rb$. Let $\{\phi_1,\cdots,\phi_m\}\subset \h_{f_0}^{0,1}\lb M_0, T^{1,0}M_0\rb$ be a basis. Then there exists a unique power series
\[
\phi(t)=\sum_{i=1}^m t_i\phi_i+\sum_{|I|\geq 2}t^I\phi_{_I}\subset A^{0,1}\lb M_0, T^{1,0}M_0\rb
\]
such that $\phi(t)$ is convergent, when $\eps$ is small, and satisfies the equations
\begin{eqnarray}\label{mc}
\begin{cases}
\bar\pa_0\phi(t)=\frac 12 \lsb \phi(t),\phi(t)\rsb;\\
\bar\pa_{f_0}^*\phi(t)=0 ; \\
\h_{f_0}\lb\phi(t)\rb=\sum_{i=1}^m t_i\phi_i \ \!.
\end{cases}
\end{eqnarray}
Furthermore, let $\XX=M_0\times B$ be the smooth manifold with the projection map $\pi:\XX\to B$,  and
\[
\Omega_{p,t}^{1,0} \ \!\XX=\pi^* T_t^{1,0}B\oplus \lb I+\phi(t)\rb\Omega_p^{1,0}M_0
\]
for any point $\lb p,t\rb\in \XX$. Then $\XX$ is a complex manifold; $\lb \XX,B,\pi\rb$ is a Kuranishi family of $M_0$;  and $t$ is a flat coordinate system, unique up to the choice of the metric $\omega_0$ in $2\pi c_1\lb M_0\rb$ and affine transformations.
\end{theorem}

In \cite{sundeform1}, the divergence gauge was introduced by the second author to study complex deformations of \ke manifolds of general type and pluricanonical forms. In particular, it was shown that, with respect to \ke metrics, the divergence gauge is equivalent to the Kuranishi gauge. Such an equivalence of the twisted Kuranishi and the divergence gauges with respect to general \ka metrics on Fano manifolds was proved by the third author \cite{zhang} in 2014; see also {\it Remark 5} in \cite{csyz2021}.

\begin{theorem}\label{eqgau}
Let $\lb M_0,\omega_0\rb$ and $f_0$ be as above. Suppose  $\phi\in A^{0,1}\lb M_0, T^{1,0}M_0\rb$ satisfies $\bar\pa_0\phi=\frac 12 \lsb\phi,\phi\rsb$.  Then
\[
\bar\pa_{f_0}^*\phi=0 \ \ \ \text{if and only if}\ \ \ \dd_{f_0}\phi=0.
\]
Furthermore, $\phi\lrcorner\omega_0=0$ when either one of these equivalent conditions is imposed.
\end{theorem}

In addition, it follows from the proof of the above theorem that

\begin{corollary}\label{ul}
Let $\lb M_0,\omega_0\rb$ and $f_0$ be as above, and let $\phi\in A^{0,1}\lb M_0, T^{1,0}M_0\rb$ such that $\bar\pa_0\phi=0$ and $\phi\lrcorner\omega_0=0$. Then $\bar\pa_0\lb \dd_{f_0}\phi\rb=0$ and
\[
\bar\pa_{f_0}^*\phi=-g^{i\bar j}\lb \dd_{f_0}\phi\rb_{\bar j}\frac{\pa}{\pa z_i} .
\]
\end{corollary}

\section{The Weil-Petersson Metric}\label{wpmc}

In this section, we first recall the \wpm metric defined on the space of Fano \ke manifolds and then explain Donaldson's construction of his new metric on $\mathcal J$ and the associated moment map. We then define the \wpm metric on the space of Fano \kr solitons and prove a necessary and sufficient condition on when such a \wpm metric is independent  of the choices of \kr soliton metrics. We also show that the \wpm metric is \ka when it defines a metric on the Kuranishi space. Finally, we establish the first and second order deformation of Fano \kr solitons and show that, essentially, the first effective term in deforming \kr solitons leads to the \wpm metric.

\subsection {The Weil-Petersson metric on the space of Fano \ke manifolds}
First of all, let us recall the \wpm metric on the Kuranishi space of Fano \ke manifolds and some key results proved in \cite{csyz2021}.
Let $\lb M_0,\omega_0\rb$ be a Fano \ke manifold and let $\lb\XX,B,\pi\rb$ be the Kuranishi family of $\lb M_0,\omega_0\rb$ as constructed in Theorem \ref{twku} (with $f_0$ being a constant potential function).  As shown in \cite{csyz2021}, there is a deep relationship between the existence of \ke metrics on each small deformation $M_t$ of $M_0$ and the stability of the action of $\au\lb M_0\rb$ on the Kuranishi space $B$.

\begin{theorem} {\bf (Cao-Sun-Yau-Zhang \cite{csyz2021})} \label{fkemain}
Let $\lb M_0,\omega_0\rb$ be a Fano \ke manifold and let $\lb \XX, B,\pi\rb$ be the Kuranishi family with respect to $\omega_0$. By shrinking $B$ if necessary, the following statements are equivalent:

\smallskip
\begin{enumerate}
\item $M_t$ admits a \ke metric for each $t\in B$;
\smallskip

\item The dimension $h^0\lb M_t, T^{1,0}M_t\rb$ of the space of holomorphic vector fields on $M_t$ is independent of $t$ for all $t\in B$;

\smallskip

\item The automorphism group $\text{Aut}_0\lb M_t\rb$ is isomorphic to $\text{Aut}_0\lb M_0\rb$ for all $t\in B$.
\end{enumerate}
In particular, the above conditions are also equivalent to the condition that the action of $\au\lb M_0\rb$ on $B$ is trivial.
\end{theorem}

The \wpm metric on the Kuranishi space $B$ is defined as

\[
\left\langle \frac{\pa}{\pa t_i},  \frac{\pa}{\pa t_j}\right\rangle_{WP}(t)=\int_{M_t} \left\langle \h_t\lb KS_t\lb  \frac{\pa}{\pa t_i}\rb\rb, \h_t\lb KS_t\lb  \frac{\pa}{\pa t_j}\rb\rb \right\rangle_{\omega_t} {\omega_t^{[n]}},
\]
where $\omega_t$ is a \ke metric on $M_t$, $KS_t$ is the Kodaira-Spencer map and $\h_t$ is the harmonic projection with respect to $\omega_t$. Since \ke metrics on $M_t$ are non-unique when $\au\lb M_t\rb$ is nontrivial, a natural question is whether the above \wpm metric depends on the choice of \ke metrics or not. In \cite{csyz2021}, we showed the following result.

\begin{theorem} {\bf (Cao-Sun-Yau-Zhang \cite{csyz2021})} \label{kewell}
The \wpm metric on $B$ is independent of the choice of \ke metrics if and only if one of the equivalent conditions in Theorem \ref{fkemain} holds.
\end{theorem}

We remark that the above results fit in with the moment map picture of Donaldson-Fujiki. Another description of the \wpm metric in the \ke case is Schumacher's work on the Deligne pairing \cite{sch}.
Let $\pi:\XX\to B$ be a family of Fano manifolds with a smooth family of \ke metrics. The volume forms of these \ke metrics define a metric on the relative anticanonical bundle $K_{\XX/B}^{-1}$. In this case, the \wpm metric is just the canonical metric on the Deligne pairing $\langle K_{\XX/B}^{-1},\cdots,K_{\XX/B}^{-1}\rangle$, up to a multiplicative constant. We note that the smoothness assumptions in Schumacher's work are in fact related to the actions of the automorphism groups of fibers on the base space $B$ as discussed in \cite{csyz2021}.

\subsection {Weil-Petersson type metric on the space of Fano \kr solitons}
In \cite{do2015}, Donaldson constructed a new \ka structure on the space of complex structures on the underlying smooth manifold of a Fano manifold which is better adapted to the \ke geometry. This work of Donaldson sheds some light on the appropriate definition of \wpm type metrics on the space of \kr solitons.  We first describe this \ka structure.

Let $\lb M, L\rb$ be the background smooth pair of a Fano manifold $M_0$ and its anti-canonical bundle $K^{-1}_{M_0}$. Let $\omega$ be a symplectic form on $M$ representing the class $2\pi c_1\lb M_0\rb$ and let $h$ be a metric on $L$ with a compatible connection $\nabla$ such that $curv\lb\nabla\rb=-\ii\omega$. Now we let
\[
\mathcal J=\lfb \text{almost complex structures on}\ M \ \text{compatible with}\ \omega\rfb
\]
and $\ji\subset\mathcal J$ be the space of integrable ones. Each element $J\in\mathcal J$ determines an $L$-valued $n$-form $\ga\in \Omega^n\lb M, L\rb$, unique up to scaling, and $J$ is integrable if and only if $d_\nabla\ga=0$. Explicitly, if $J\in\ji$ then, up to a multiplicative constant, $\ga$ coincides with $dz_1\wedge\cdots dz_n\otimes \lb \frac{\pa}{\pa z_1}\wedge\cdots  \frac{\pa}{\pa z_n}\rb$  for any local holomorphic coordinates $(z_1,\cdots,z_n)$ on the complex manifold $M_J=\lb M, J\rb$. We let $\hat{\mathcal J}^{int}\subset \Omega^n\lb M, L\rb$ be the space of all such $\ga$ and denote by $T_\ga$ the tangent space of $\hat{\mathcal J}^{int}$ at $\ga$. Naturally $\hat{\mathcal J}^{int}$ is a $\C^*$-bundle over $\ji$. The metric $h$ on $L$ leads to a natural $L^2$ metric on $\Omega^n\lb M,L\rb$. For $\ga,\gb\in \Omega^n\lb M,L\rb$, set
\begin{eqnarray}\label{dowp10}
\langle\langle \ga, \gb\rangle\rangle=-c_n\int_M \lb \ga\wedge\bar\gb\rb_h
\end{eqnarray}
where $c_n=\ii$ if $n$ is odd, and $c_n=1$ if $n$ is even. Donaldson showed that, for $\ga\in \Omega^n\lb M, L\rb$ with $d_\nabla\ga=0$, $\langle\langle\cdot,\cdot\rangle\rangle$ is positive definite on the orthogonal complement of $\ga$ in $T_\ga$ and descends to a \ka metric on $\ji$.

Now we can rewrite Donaldson's metric in terms of Beltrami differentials (see also \cite{prsa21}). For each $J\in\ji$, we have the identification
\[
T_J^{1,0}\ji\cong \lfb \phi\in A^{0,1}\lb M_J, T^{1,0}M_J\rb\big | \ \bar\pa_J\phi=0,\ \ \phi\lrcorner\omega=0\rfb.
\]
Let $\Omega_J$ be the unique volume form on $M$ with
\[
\rc_J\lb \Omega_J\rb=\omega\ \ \ \text{and}\ \ \ \int_M \Omega_J=\int_M {\omega^{[n]}}.
\]
Namely, if we let $f_J$ be the normalized Ricci potential of the \ka manifold $\lb M_J,\omega\rb$ then $\Omega_J=e^{f_J}{\omega^{[n]}}$.

By Corollary \ref{ul}, we know that $\bar\pa_J\dd_{f_J}\phi=0$ for each $\phi\in T_J^{1,0}\ji$. Since $h^{0,1}\lb M_J\rb=0$, there exists a unique function $\xi_\phi\in C^\infty\lb M, \C\rb$ such that
\[
\dd_{f_J}\phi=\bar\pa_J\xi_\phi\ \ \ \text{and}\ \ \ \int_{M_J}\xi_\phi\Omega_J=0.
\]
Then, for any $\phi,\psi\in T_J^{1,0}\ji$, the Donaldson metric is given by
\begin{eqnarray}\label{dowp20}
\langle\langle\phi,\psi\rangle\rangle =\int_{M_J}\lb \text{Tr}\lb\phi\bar\psi\rb-\xi_\phi\bar{\xi_\psi}\rb\Omega_J.
\end{eqnarray}
We note that if $g_J$ is the \ka metric on $M_J$ with \ka form $\omega$ then
\[
 \text{Tr}\lb\phi\bar\psi\rb=\langle \phi,\psi\rangle_{g_J}
\]
because $\phi\lrcorner\omega=\psi\lrcorner\omega=0$. Now we let $G=\text{Symp}_0\lb M,\omega\rb$. Since $M_0$ is a Fano manifold, we know that $H^1\lb M\rb=0$. It follows that $G$ consists of Hamiltonian diffeomorphisms. We identify the Lie algebra of $G$ by
\[
\g\cong \lfb u\in C^\infty\lb M,\R\rb\bigg |\ \int_M e^u{\omega^{[n]}}=\int_M{\omega^{[n]}}\rfb.
\]
Then the dual of $\g$ can be identified with
\[
\g^*\cong\lfb \Omega\in \Omega^{2n}\lb M,\R\rb\bigg | \ \int_M \Omega=\int_M {\omega^{[n]}}\rfb.
\]
In this case, the moment map $\mu$ of the action of $G$ on $\lb \ji, \langle\langle\cdot,\cdot\rangle\rangle\rb$ is given by $\mu\lb J\rb=\Omega_J$.

In view of Donaldson's construction above and the twisted Kuranishi-divergence gauge discussed in Section 2, it is natural to define the \wpm metric on the deformation space of Fano \kr solitons to be the restriction of the Donaldson metric to a Kuranishi space with respect to a twisted Kuranishi-divergence gauge determined by \kr solitons.

Precisely, let $p:\mathcal Y\to C$ be a holomorphic family of Fano manifolds such that each fiber $Y_s=p^{-1}(s)$ admits a \kr soliton metric $\omega_s$ and the family $\lfb\omega_s\rfb$ depends smoothly on $s$. Then the  \wpm metric is the Hermitian metric on $C$ defined in the following way.

\medskip
% \begin{definition}
\noindent{\bf Definition 3.1.} \label{wpdef}
For any $s\in C$ and any tangent vectors $u,v\in T_s^{1,0}C$, consider
\[
\text{KS}_s\lb u\rb, \text{KS}_s\lb v\rb\in H^{0,1}\lb Y_s, T^{1,0}Y_s\rb
\]
where KS is the Kodaira-Spencer map. Let $f_s$ be the (normalized) Ricci potential of $\omega_s$ and let $\phi,\psi\in \h_{f_s}^{0,1}\lb Y_s, T^{1,0}Y_s\rb$ be the $f_s$-harmonic representatives of
$\text{KS}_s\lb u\rb$ and $\text{KS}_s\lb v\rb$, respectively. Then the {\it \wpm inner product} of $u$ and $v$ is
\begin{eqnarray}\label{wpdefn10}
\langle u,v\rangle_{_{WP}}:=\int_{Y_s}\text{Tr}\lb\phi\bar\psi\rb\ e^{f_s} {\omega_s^{[n]}}.
\end{eqnarray}
% \end{definition}
\medskip

\noindent {\it Remark 3.1.} \label{dep}
In the above definition, if the family $\mathcal Y$ is a Kuranishi family of $Y_s$ then the \wpm metric at $s\in C$ is essentially a Hermitian metric on the cohomology $H^{0,1}\lb Y_s, T^{1,0}Y_s\rb$. Note that in general it may depend on the choice of the \kr soliton metric on $Y_s$.
% \end{remark}

\medskip
Now we investigate the relation between the Donaldson metric and the \wpm metric for general Fano manifolds.

\begin{proposition}\label{dowp}
Let $\lb M,\omega\rb$ be a Fano manifold such that $\rc\lb\omega\rb=\omega+\ii\pa\bar\pa f$. For any (cohomology classes) $A,B\in H^{0,1}\lb M, T^{1,0}M\rb$ and $\phi,\psi\in A^{0,1}\lb M, T^{1,0}M\rb$, with $\bar\pa\phi=\bar\pa\psi=0$, $\phi\lrcorner\omega=\psi\lrcorner\omega=0$ and $\lsb\phi\rsb=A$, $\lsb\psi\rsb=B$, we have
\[
\langle A,B\rangle_{_{WP}}=\langle\langle\phi,\psi\rangle\rangle+\int_M\lb  \lb 1+\Delta_f\rb^{-1}\xi_\phi\rb\bar{\xi_\psi}\ e^f {\omega^{[n]}},
\]
where $\langle\langle\phi,\psi\rangle\rangle$ is the Donaldson metric as defined in equations \eqref{dowp10} and \eqref{dowp20}.
\end{proposition}

\begin{proof}

Since $\bar\pa\phi=0$, by Hodge theory, we know that there exists a unique smooth vector field $W_\phi\in A^0\lb M, T^{1,0}M\rb$ such that
\[
\begin{cases}
\phi-\bar\pa W_\phi=\h_f\lb \phi\rb\\
W_\phi\bot_{L^2} H^0\lb M, T^{1,0}M\rb.
\end{cases}
\]
By Theorems \ref{twku} and \ref{eqgau}, we know that $\dd_f\h_f\lb\phi\rb=0$. It follows that
\[
\dd_f\phi-\dd_f\bar\pa W_\phi=\dd_f\h_f\lb\phi\rb=0
\]
and
\[
\dd_f\bar\pa W_\phi=\bar\pa\dd_f W_\phi-\ii W_\phi\lrcorner\omega.
\]
By Corollary \ref{ul}, we know that $\bar\pa\lb\dd_f\phi\rb=0$ and thus $\bar\pa\lb W_\phi\lrcorner\omega\rb=0$, which implies that there exists a unique smooth function $u_\phi$ such that
\[
\begin{cases}
W_\phi\lrcorner\omega=\ii \bar\pa u_\phi\\
\int_M u_\phi e^f{\omega^{[n]}}=0.
\end{cases}
\]
It then follows that $W_\phi=\nabla^{1,0}u_\phi$, $\dd_f W_\phi=\Delta_f u_\phi$, and
\[
\bar\pa\xi_\phi=\dd_f\phi=\bar\pa\lb\dd_f W_\phi\rb+\bar\pa u_\phi=\bar\pa \lb \lb 1+\Delta_f\rb u_\phi\rb.
\]
By the normalization of $\xi_\phi$ and $u_\phi$, we conclude that $\xi_\phi= \lb 1+\Delta_f\rb u_\phi$.

Similarly, we have $\xi_\psi= \lb 1+\Delta_f\rb u_\psi$. Thus,
\begin{align*}
\begin{split}
\langle A,B\rangle_{_{WP}}=& \langle\langle \h_f\lb\phi\rb,\h_f\lb\psi\rb  \rangle\rangle\\
=&\int_M \text{Tr}\lb\phi\bar\psi\rb  \ e^f{\omega^{[n]}} -
\int_M \langle\bar\pa W_\phi,\bar\pa W_\psi\rangle_\omega\  e^f{\omega^{[n]}}\\
=& \int_M \text{Tr}\lb\phi\bar\psi\rb\  e^f{\omega^{[n]}} -\int_M \Delta_f u_\phi \bar{\lb 1+\Delta_f\rb u_\psi}\ e^f{\omega^{[n]}}\\
=& \langle\langle \phi,\psi \rangle\rangle+\int_M u_\phi\bar{\xi_\psi}\ e^f{\omega^{[n]}}\\
=&\langle\langle\phi,\psi\rangle\rangle+\int_M\lb  \lb 1+\Delta_f\rb^{-1}\xi_\phi\rb\bar{\xi_\psi}\ e^f{\omega^{[n]}}.
\end{split}
\end{align*}
\end{proof}

For a Fano manifold $M$ admitting a \kr soliton metric, we know that \kr soliton metrics on $M$ are not unique. However, for any two such soliton metrics $\omega_1$ and $\omega_2$ on $M$, by a result of Tian-Zhu \cite{tianzhu}, there exists a biholomorphism $\sigma\in\au\lb M\rb$ such that $\omega_2=\sigma^*\omega_1$. In order to obtain a metric on the moduli space of \kr solitons over $M$, it is  natural to require that the \wpm metric defined by (\ref{wpdefn10}) is independent of the choices of \kr solitons in the spirit of Theorem \ref{kewell}. If so, we shall say that the \wpm metric is well-defined.

Our following result provides a necessary and sufficient condition on the well-definedness of the \wpm metric.

\begin{theorem}\label{wd}
Let $M$ be a Fano manifold that admits a \kr soliton metric. Then the \wpm metric is well-defined if and only if the action of $\au\lb M\rb$ on the Kuranishi space of $M$ is trivial.
\end{theorem}

\begin{proof}

Let $\omega_0$ be any \kr soliton metric on $M$ with normalized Ricci potential $f_0$, and let $\Omega_0=e^{f_0}{\omega_0^{[n]}}$ be the weighted volume form. For any smooth family $\lfb\sigma_s\rfb_{ -\eps<s<\eps}\subset\au\lb M\rb$ with $\sigma_0=Id$, let $\omega_s=\sigma_s^*\omega_0$. Then $\omega_s$ is a \kr soliton metric with normalized Ricci potential $f_s$. We denote $\Omega_s=e^{f_s}{\omega_s^{[n]}}$.

For any cohomology classes $A,B\in H^{0,1}\lb M, T^{1,0}M\rb$, we let $\phi_s,\psi_s\in \h_{f_s}^{0,1}\lb M, T^{1,0}M\rb$ with $\lsb\phi_s\rsb=A$ and $\lsb\psi_s\rsb=B$.
Since $\omega_0$ is an arbitrary \kr soliton metric, the \wpm metric is well-defined if and only if
\[
\frac{d}{ds}\bigg |_{s=0}\int_M \tr \lb\phi_s\bar\psi_s\rb\Omega_s=0.
\]
Let $w=\frac{d}{ds}\bigg |_{s=0}\sigma_s$. Then  $w$ is a real holomorphic vector field on $M$. Since $\lsb\phi_s\rsb=\lsb\phi_0\rsb$ and $\lsb\psi_s\rsb=\lsb\psi_0\rsb$,  there exist smooth families of vector fields $X(s),Y(s)\subset A^0\lb M, T^{1,0}M\rb$, respectively,  with $X(0)=Y(0)=0$, such that
 \[\phi_s=\phi_0+\bar\pa X(s) \quad {\text and} \quad  \psi_s=\psi_0+\bar\pa Y(s).\]
It follows from a direct computation that

\begin{align*}
\begin{split}
\frac{d}{ds}\bigg |_{s=0}\int_M \tr \lb\phi_s\bar\psi_s\rb\Omega_s
=&\int_M \left \langle \bar\pa \lb \frac{d}{ds}\bigg |_{s=0} X(s)\rb, \psi_0\right\rangle_{\omega_0}\Omega_0\\
&+\int_M \left \langle \phi_0, \bar\pa \lb \frac{d}{ds}\bigg |_{s=0} Y(s)\rb\right\rangle_{\omega_0}\Omega_0\\
&+\int_M \tr\lb\phi_0\bar\psi_0\rb \frac{d}{ds}\bigg |_{s=0}\\
=& \int_M \tr\lb\phi_0\bar\psi_0\rb\lb\dd_{f_0}w\rb\Omega_0.
\end{split}
\end{align*}
Replacing $w$ by $Jw$ and combining, we see that the \wpm metric is well-defined provided, for any \kr soliton metric $\omega_0$, any basis $\{\phi_1,\cdots,\phi_m\}$ of $\h_{f_0}^{0,1}\lb M, T^{1,0}M\rb$ and any holomorphic vector field $v\in H^0\lb M, T^{1,0}M\rb$, we have
\begin{eqnarray}\label{ip05}
\int_M \lb \dd_{f_0}v\rb\tr\lb\phi_i\bar\phi_j\rb\Omega_0=0,  \ \ \   \ 1\leq i, j\leq m.
\end{eqnarray}
Namely,
\begin{eqnarray}\label{ip10}
\lfb \tr\lb\phi_i\bar\phi_j\rb\mid 1\leq i,j\leq m\rfb \bot_{L^2}\Lambda_{f_0}^1 ,
\end{eqnarray}
where we use the volume form $\Omega_0$ to define the $L^2$ inner product. By integration  by parts,  we see that equation \eqref{ip05} is equivalent to
\begin{eqnarray}\label{ip20}
\int_M  \left\langle \lsb v,\phi_i\rsb, \phi_j\right\rangle_{\omega_0}\Omega_0=0.
\end{eqnarray}
Since $\lsb v,\phi_i\rsb=L_v\phi_i$, $\bar\pa\lb L_v\phi_i\rb=0$ and $\phi_j$ is $f_0$-harmonic, we see that equation \eqref{ip20} is in turn equivalent to the vanishing of the cohomology classes
\begin{eqnarray}\label{ip30}
\lsb L_v\phi_i\rsb=0 , \ \ \ \text{for}\  \ 1\leq i\leq m.
\end{eqnarray}
Now we let $\pi:\XX\to B$ be a Kuranishi family of $M$ such that $\pi^{-1}(0)\cong M$. Then $\au\lb M\rb$ acts on $B$ and fixes $0$. Hence $\au\lb M\rb$ acts on $T_0^{1,0}B\cong H^{0,1}\lb M, T^{1,0}M\rb$ linearly, and the corresponding action of its Lie algebra, $\text{Lie}\lb \au\lb M\rb\rb\cong H^0\lb M, T^{1,0}M\rb$, on $H^{0,1}\lb M, T^{1,0}M\rb$ is given by $v\lb \lsb\phi\rsb\rb=\lsb L_v\phi\rsb$. Thus equation \eqref{ip30} is equivalent to the condition that the action of $\au\lb M\rb$ on $B$ is trivial.

\end{proof}

A direct corollary, as in \cite{csyz2021},   is the following sufficient condition on the well-definedness of the \wpm metric.

\begin{corollary}\label{nj}
Let $\lb \XX,B,\pi\rb$ be a Kuranishi family of $M\cong M_0=\pi^{-1}(0)$ such that $M$ admits a \kr soliton. If $h^0\lb M_t, T^{1,0}M_t\rb$ remains constant as $t$ varies in $B$ then the \wpm metric is well-defined.
\end{corollary}

\medskip
\noindent {\it Remark 3.2.\ }
% \begin{remark}
As stated in Theorem \ref{fkemain}, the conditions of the trivial action of $\au\lb M\rb$ on the Kuranishi space of $M$ and the constancy of $h^0\lb M_t, T^{1,0}M_t\rb$ in $t$ are equivalent in the case of Fano \ke metrics, and they are also equivalent to the existence of \ke metrics on each fiber $M_t$. Such equivalence relations in the  case of Fano \kr solitons will be discussed in a forthcoming paper.

The automorphism group of a shrinking \kr soliton is not reductive in general. However, it is not hard to show that if the action of the group generated by the soliton vector field on the Kuranishi space is trivial, then the action of the unipotent radical is also trivial. This is clear in the level of Lie algebra and is followed from Mabuchi's work on the Calabi-Matsushima type theorem \cite{macama} and the Jacobi identity. 
% \end{remark}

\medskip
Now we assume that there exists a smooth family of \kr solitons on a Kuranishi family $\lb \XX, B,\pi\rb$ of Fano manifolds and the \wpm metric is well-defined at each point of $B$. In this case, the \wpm metric is a Hermitian metric on the Kuranishi space $B$.

In general, a basic question is whether  \wpm metrics are K\"ahler. In the case of Riemann surfaces, the K\"ahlerian property of the \wpm metric was stated by Weil \cite{weil} and proved by Ahlfors \cite{ahlforsa}. In \cite{koiso83},  Koiso showed that the \wpm metrics on the spaces of \ke manifolds with nonzero first Chern class are K\"ahler.
Later,  Siu \cite{siu} provided a simpler proof by using the canonical lift. In \cite{nan}, Nannicini showed that the \wpm metrics on the moduli spaces of polarized Calabi-Yau manifolds are K\"ahler. In \cite{sch}, Schumacher gave a  unified proof by using the harmonic lift which depends on the variation of \ke metrics.

In the case of Fano \kr solitons, to compute the derivatives of the \wpm metric and check its K\"ahlerian property, we need to study the deformation of \kr solitons. When the \wpm metric is well-defined, we have freedom to impose any gauge and to modify the given family of \kr solitons.

Now let $\lb \XX, B,\pi\rb$ be a Kuranishi family of $M_0=\pi^{-1}(0)$ as constructed in Theorem \ref{twku} with $\phi(t)$ the unique solution to \eqref{mc} and let $\lfb\omega_t\rfb_{t\in B}$ be a smooth family of \ka metrics such that $\omega_t$ is a \kr soliton metric on $M_t$ with normalized Ricci potential $f_t$. We impose the $f_0$-twisted Kuranishi-divergence gauge on $\XX$.
As before, for each $t\in B$, we let $\Omega_t=e^{f_t}{\omega_t^{[n]}}$ and also $X_t=\nabla_t^{1,0}f_t\in H^0\lb M_t, T^{1,0}M_t\rb$ be the soliton vector field on $M_t$. It follows that $\text{Im}\lb X_t\rb$ is a Killing field on $\lb M_t,\omega_t\rb$. We denote by $T$ the subgroup of $\text{Isom}_0\lb M_0,\omega_0\rb$ generated by $\text{Im}\lb X_0\rb$.

The twisted Kuranishi-divergence gauge allows us to identify each $M_t$ with $M_0$ as smooth manifolds and we can view $\lfb\Omega_t\rfb_{t\in B}$ as a family of volume forms on $M_0$ and thus there exists a smooth function $\rho=\rho\lb t,z\rb\in C^\infty\lb B\times M_0,\ \R\rb$ such that
\begin{eqnarray}\label{ip35}
\Omega_t=e^\rho\det\lb I-\phi(t)\bar{\phi(t)}\rb\Omega_0.
\end{eqnarray}
Here, we keep the factor $\det\lb I-\phi(t)\bar{\phi(t)}\rb$ in (\ref{ip35}) simply for the convenience of computations. It follows that $\rho(0,z)=0$ and we have
\begin{eqnarray}\label{ip40}
\rho=\sum_{i=1}^m t_i\mu_i+\sum_{j=1}^m \bar t_j\bar\mu_j+O\lb |t|^2\rb,
\end{eqnarray}
where $m=h^{0,1}\lb M_0, T^{1,0}M_0\rb$ and $\mu_i\in C^\infty\lb M_0,\C\rb$, $1\le i\le m$,  are smooth functions. In this situation, we have

\begin{theorem}\label{1st}
The Kuranishi family $\lb\XX, B,\pi\rb$ is $T$-equivariant and the potential functions $f_t$ have the expansion
\begin{eqnarray}\label{ip45}
f_t=f_0+\sum_{i=1}^m t_i\lb 1+\Delta_0\rb\mu_i+\sum_{j=1}^m \bar t_j\lb 1+\Delta_0\rb\bar\mu_j+O\lb |t|^2\rb.
\end{eqnarray}
Furthermore, each $\mu_i$ satisfies the equations
\begin{eqnarray}\label{ip50}
\begin{cases}
\lb 1+\Delta_{f_0}\rb\lb 1+\bar\Delta_{f_0}\rb \mu_i=0 ,\\
\lb 1+\Delta_{f_0}\rb\lb 1+\bar\Delta_{f_0}\rb \bar\mu_i=0 .
\end{cases}
\end{eqnarray}
\end{theorem}

\begin{proof}

By Theorems \ref{twku} and \ref{eqgau}, we know that $\dd_{f_0}\phi(t)=0$ and $\phi(t)\lrcorner\omega_0=0$ for each $t$. By using these results and the facts that  $\omega_t=-\ii\pa_t\bar\pa_t\log \Omega_t$ and
$f_t=\log \lb{\Omega_t}/{\omega_t^{[n]}}\rb$, formula \eqref{ip45} follows from a direct computation.

Now, since $\omega_t$ is a \kr soliton metric and thus $\bar\pa_t\nabla_t^{1,0}f_t=0$ for each $t$, we have
\[\frac{\pa}{\pa t_i}\bigg |_{t=0}\bar\pa_t\nabla_t^{1,0}f_t=0 \quad  \text {and} \quad  \frac{\pa}{\pa \bar t_j}\bigg |_{t=0}\bar\pa_t\nabla_t^{1,0}f_t=0, \quad  \text {for each} \  i, j.\]
Then the first equation leads to
\begin{eqnarray}\label{ip60}
L_{\text{Im}\lb X_0\rb}\phi_i +\bar\pa_0\nabla_0^{1,0}\lb 1+\bar\Delta_{f_0}\rb\mu_i=0
\end{eqnarray}
for each $i$, and the second equation leads to
\begin{eqnarray}\label{ip70}
\bar\pa_0\nabla_0^{1,0}\lb 1+\bar\Delta_{f_0}\rb\bar\mu_j=0
\end{eqnarray}
for each $j$.

By Theorem \ref{eqgau}, we know that $\phi_i\lrcorner\omega_0=0$ and $\dd_{f_0}\phi_i=0$. A direct computation then shows that
\[
L_{\text{Im}\lb X_0\rb}\phi_i\in\h_{f_0}^{1,0}\lb M_0, T^{1,0}M_0\rb
\]
is $f_0$-harmonic. Thus, by equation \eqref{ip60}, we know that
\begin{eqnarray}\label{ip80}
L_{\text{Im}\lb X_0\rb}\phi_i =0
\end{eqnarray}
for each $i$. This implies that the action of $T$ on $\h_{f_0}^{1,0}\lb M_0, T^{1,0}M_0\rb$ is trivial. Since $T\subset \text{Isom}_0\lb M_0,\omega_0\rb$, by the uniqueness part of Theorem \ref{twku} we know that $T$ preserves $\phi(t)$ for each $t\in B$ and thus $T\subset\au\lb M_0\rb$.

Finally, formula \eqref{ip50} follows from equations \eqref{ip60}-\eqref{ip80} and Theorem \ref{e1}.

\end{proof}

% \begin{remark}
\noindent{\it Remark 3.3.} \label{eqfamily}
The $T$-equivariance of the Kuranishi family allows us to modify the family $\lfb\omega_t\rfb$ of \kr solitons. In fact, by averaging the family $\lfb\omega_t\rfb$ over $T$ with respect to the Haar measure, we get a $T$-invariant family $\lfb\wt\omega_t\rfb$ of \kr soliton metrics on $\XX$ which restrict to $\omega_0$ on $M_0$ since $\omega_0$ is $T$-invariant.
% \end{remark}
\medskip

From now on we will assume the family $\lfb\omega_t\rfb$ of \kr soliton metrics is $T$-invariant. In this case, the functions $\mu_i$'s in the expansion \eqref{ip40} are very special, hence $\rho$ takes a much simpler form.

\begin{corollary}\label{msi}
Let $\lfb\omega_t\rfb$ be a $T$-invariant family of \kr soliton metrics on the Kuranishi family $\lb\XX,B,\pi\rb$. Then, each $\mu_i$ in the expansion \eqref{ip40} has the property that
\[
\mu_i,\ \bar\mu_i\in\Lambda_{f_0}^1.
\]
\end{corollary}

\begin{proof}

Since each $\omega_t$ is $T$-invariant, we know that $\rho$ is $T$-invariant and thus $\text{Im}\lb X_0\rb\lb\rho\rb=0$. This implies $\text{Im}\lb X_0\rb\lb\mu_i\rb=0$ and $\text{Im}\lb X_0\rb\lb\bar\mu_j\rb=0$. It then follows that

\begin{align*}
\begin{split}
\lb 1+\bar\Delta_{f_0}\rb \mu_i= & \lb 1+\Delta_{0}\rb \mu_i+X_0\lb\mu_i\rb \\
=& \lb 1+\Delta_{0}\rb \mu_i+\bar{X_0}\lb\mu_i\rb \\
= & \lb 1+\Delta_{f_0}\rb \mu_i,
\end {split}
\end{align*}
hence
\begin{align*}
\begin{split}
0=& \int_{M_0} \lb\lb 1+\Delta_{f_0}\rb\lb 1+\bar\Delta_{f_0}\rb \mu_i\rb\bar{\mu_i}\ \Omega_0 \\
= & \int_{M_0}\lb 1+\Delta_{f_0}\rb\mu_i\bar{\lb 1+\Delta_{f_0}\rb\mu_i}\ \Omega_0.
\end {split}
\end{align*}
Thus $\lb 1+\Delta_{f_0}\rb\mu_i=0$, and the same argument works for $\bar\mu_i$.

\end{proof}

%\begin{remark}
\noindent {\it Remark 3.4.} \label{reduc} \
By the above corollary, we know that $\text{Re}\lb\mu_i\rb \in\Lambda_{f_0}^1$ and $ \text{Im}\lb\mu_i\rb \in\Lambda_{f_0}^1$.  It follows that, for each $i$, the holomorphic vector field $\nabla_0^{1,0}\mu_i$ lies in the reductive part of $H^0\lb M_0, T^{1,0}M_0\rb$ that is determined by $\omega_0$. In particular, $\nabla_0^{1,0}\mu_i$ commutes with $X_0$ and $\text{Im}\lb X_0\rb$ since $X_0$ lies in the center of the reductive part.

% \end{remark}
\medskip

By using Theorem \ref{1st}, we can now prove the K\"ahlerian property of the \wpm metric.

\begin{theorem}\label{wpka}
Assume that each small deformation of a Fano \kr soliton  admits a \kr soliton metric and the \wpm metric is well-defined. Then the \wpm metric on the Kuranishi space is K\"ahler. Furthermore, the flat coordinate system described in Theorem \ref{twku} is a normal coordinate system of the \wpm metric at $0$.
\end{theorem}

\begin{proof}

We fix a \kr soliton $\lb M_0, \omega_0\rb$ with normalized Ricci potential $f_0$ and let $\lb \XX,B,\pi\rb$ be the Kuranishi family of $M_0$ constructed in Theorem \ref{twku}. By the assumption in Theorem~\ref{wpka} and Remark \ref{eqfamily}, we can extend $\omega_0$ to a $T$-invariant family of soliton metrics $\lfb\omega_t\rfb$ on $\XX$ such that $\omega_t$ is a \kr soliton on $M_t$ for each $t\in B$. Let
\[
h_{i\bar j}(t)=\left\langle \frac{\pa}{\pa t_i}, \frac{\pa}{\pa t_j}\right\rangle_{_{WP}}(t)\equiv \int_{M_t}\text{Tr} \lb \h_{f_t} \lb KS_t\lb \frac{\pa}{\pa t_i}\rb\rb
\bar{\h_{f_t} \lb KS_t\lb \frac{\pa}{\pa t_j}\rb\rb }\rb\Omega_t.
\]
To prove the theorem, it is enough to show that $\frac{\pa}{\pa t_k}h_{i\bar j}(0)=0$ for all $i,j,k$.

We first consider the Kodaira-Spencer class $KS_t\lb \frac{\pa}{\pa t_i}\rb$. Let
\[
\psi_i\lb t\rb=\lb  \frac{\pa}{\pa t_i}\phi(t)\rb \lb I-\bar{\phi(t)}\phi(t)\rb^{-1}.
\]
Then a direct computation shows that $\bar\pa_t\psi_i(t)=0$ and $KS_t\lb \frac{\pa}{\pa t_i}\rb=\lsb \psi_i(t)\rsb$. To find the harmonic representative $\h_{f_t}\lb \psi_i(t)\rb$, we recall the linear isomorphism
\[
\tau_t:A^0 \lb M_0, T^{1,0}M_0\rb\to A^0 \lb M_t, T^{1,0}M_t\rb
\]
that was defined in \cite{csyz2021},
\[
\tau_t(v)=\lb I-\phi(t)\bar{\phi(t)}\rb^{-1}(v)-\bar{\phi(t)}\lb \lb I-\phi(t)\bar{\phi(t)}\rb^{-1}(v)\rb.
\]
By using this isomorphism, there exists a smooth family $\lfb V_i(t)\rfb_{t\in B}\subset A^0 \lb M_0, T^{1,0}M_0\rb$ such that $V_i(0)=0$ and
\[
\psi_i(t)+\bar\pa_t \lb \tau_t\lb V_i(t)\rb\rb=\h_{f_t}\lb \psi_i(t)\rb.
\]
By direct computations, we have
\begin{align*}
\begin{split}
\frac{\pa}{\pa t_k}\bigg |_{t=0}\psi_i(t)=& \frac{\pa^2}{\pa t_i\pa t_k}\bigg |_{t=0}\phi(t),\\
\frac{\pa}{\pa\bar t_l}\bigg |_{t=0}\psi_i(t)=& 0,\\
\frac{\pa}{\pa t_k}\bigg |_{t=0}\bar\pa_t \lb \tau_t\lb V_i(t)\rb\rb=& \bar\pa_0\lb \frac{\pa}{\pa t_k}\bigg |_{t=0}V_i(t)\rb,
\\
\frac{\pa}{\pa\bar t_l}\bigg |_{t=0}\bar\pa_t \lb \tau_t\lb V_i(t)\rb\rb=& \bar\pa_0\lb \frac{\pa}{\pa\bar t_l}\bigg |_{t=0}V_i(t)\rb.
\end{split}
\end{align*}
We note that if  $\nu_1,\nu_2\in A^{0,1}\lb M_t, T^{1,0}M_t\rb$ are two Beltrami differentials such that either $\nu_1\lrcorner\omega_t=0$ or  $\nu_2\lrcorner\omega_t=0$, then $\tr\lb\nu_1\bar\nu_2\rb=\langle \nu_1,\nu_2\rangle_{\omega_t}$. It follows that

\begin{align*}
\begin{split}
\frac{\pa}{\pa t_k}h_{i\bar j}(0)=& \int_{M_0} \tr\lb \frac{\pa}{\pa t_k}\bigg |_{t=0}\lb \psi_i(t)+\bar\pa_t \lb \tau_t\lb V_i(t)\rb\rb\rb\bar{\phi_j}\rb\Omega_0\\
&+  \int_{M_0} \tr \lb \phi_i \bar{ \frac{\pa}{\pa\bar t_k}\bigg |_{t=0}\lb \psi_j(t)+\bar\pa_t \lb \tau_t\lb V_j(t)\rb\rb\rb}\rb\Omega_0\\
&+ \int_{M_0} \tr \lb \phi_i \bar\phi_j\rb \lb \frac{\pa}{\pa t_k}\bigg |_{t=0}\Omega_t\rb\\
=& \int_{M_0} \left \langle \frac{\pa^2}{\pa t_i\pa t_k}\bigg |_{t=0}\phi(t),\phi_j \right \rangle_{\omega_0}\ \Omega_0
+ \int_{M_0} \left \langle\bar\pa_0\lb \frac{\pa}{\pa t_k}\bigg |_{t=0}V_i(t)\rb,\phi_j \right \rangle_{\omega_0}\Omega_0\\
&+  \int_{M_0} \left \langle\phi_i, \bar\pa_0\lb \frac{\pa}{\pa\bar t_l}\bigg |_{t=0}V_i(t)\rb \right \rangle_{\omega_0}\Omega_0
+ \int_{M_0} \tr \lb \phi_i \bar\phi_j\rb\mu_k\ \Omega_0\\
=&  \int_{M_0} \left \langle \frac{\pa^2}{\pa t_i\pa t_k}\bigg |_{t=0}\phi(t),\phi_j \right \rangle_{\omega_0}\ \Omega_0
+ \int_{M_0} \tr \lb \phi_i \bar\phi_j\rb\mu_k\ \Omega_0.
\end{split}
\end{align*}
By the construction of the power series $\phi(t)$ in Theorem \ref{twku}, we know that $\frac{\pa^2}{\pa t_i\pa t_k}\bigg |_{t=0}\phi(t)$ is in the image of $\bar\pa_{f_0}^*$ and thus the first term of the right side of the above formula vanishes. To handle the second term, by Corollary \ref{msi}, we know that $\mu_k\in\Lambda_{f_0}^1$. Since the \wpm metric is well-defined, equation \eqref{ip10} implies that the second term vanishes.

Thus the \wpm metric is \ka and the flat coordinates $(t_1,\cdots,t_m)$ is a normal coordinate system of the \wpm metric at $0$.

\end{proof}

Finally, we show that the \wpm metric is closely related to the second order deformation of \kr solitons. To see this clearly and to simplify computations, we need to modify the $T$-invariant family $\lfb\omega_t\rfb$ of \kr soliton metrics on the Kuranishi family $\XX$ once more.

We consider the holomorphic vector field $\nabla_0^{1,0}\mu_k$ for $1\le k\le m$. By Remark \ref{reduc}, we know that it is $T$-invariant. Now we extend $\nabla_0^{1,0}\mu_k$ to a holomorphic section $Z_k\in H^0\lb T^{1,0}_{\XX/B}\rb$. By the work of Kodaira, this can be done provided the dimension of the reductive part of $\au\lb M_t\rb$ remains constant as $t$ varies in $B$, and this is closely related to the existence of \kr soliton metrics on each $M_t$ in the spirit of Theorem \ref{fkemain}. Again, by the averaging trick  which leaves $\nabla_0^{1,0}\mu_k$ intact, we can assume that each $Z_k$ is $T$-invariant. Now we let $\sigma$ be the time $1$ flow of the vector field $Z=\sum_{k=1}^m t_k Z_k$ and we replace the family $\lfb \omega_t\rfb$ of \kr soliton metrics by the family $\lfb \sigma^*\omega_t\rfb$. A direct computation shows that all $\mu_k$'s as in the expansion \eqref{ip40} with respect to the new family $\lfb \sigma^*\omega_t\rfb$ vanish. We call such a family of \kr solitons a normalized family. In this case,  the function $\rho$ defined by equation \eqref{ip35} is special that it satisfies $\rho=O\lb |t|^2\rb$.

\begin{theorem}\label{2nd}
Let $\lb \XX,B,\pi\rb$ be the Kuranishi family of a \kr soliton $\lb M_0,\omega_0\rb$ given by Theorem \ref{twku}, and let $\lfb \omega_t\rfb_{t\in B}$ be a $T$-invariant normalized family of \kr solitons on $\XX$. Let $\rho$ be the deformation function defined by equation \eqref{ip35} and set
\[
\eta_{i\bar j}=\frac{\pa^2}{\pa t_i\pa\bar t_j}\bigg |_{t=0}\rho.
\]
Then we have
\begin{enumerate}

\item [(a)]
\begin{eqnarray}\label{ip90}
\frac{\pa^2}{\pa t_i\pa\bar t_j}\bigg |_{t=0} f_t=\lb 1+\Delta_{0}\rb\eta_{i\bar j}-\tr\lb\phi_i\bar\phi_j\rb,
\end{eqnarray}
where $f_t$ is the normalized Ricci potential of $\omega_t$.

\item  [(b)]
\begin{eqnarray}\label{ip100}
\lb 1+\Delta_{f_0}\rb\eta_{i\bar j}-\tr\lb\phi_i\bar\phi_j\rb \in \Lambda_{f_0}^1.
\end{eqnarray}
\end{enumerate}

In particular,
\begin{eqnarray}\label{ip110}
\int_{M_0}\eta_{i\bar j}\ \Omega_0=\left\langle \frac{\pa}{\pa t_i}, \frac{\pa}{\pa t_j}\right\rangle_{WP}(0).
\end{eqnarray}
\end{theorem}

\medskip

\begin{proof}
The proof of formula \eqref{ip90} is similar to the first order expansion of $f_t$ in the proof of Theorem \ref{1st}.
To prove \eqref{ip100}, by the \kr soliton equation, we note that
\[\frac{\pa^2}{\pa t_i\pa\bar t_j}\bigg |_{t=0}\bar\pa_t\nabla_t^{1,0}f_t=0.\] A simple computation shows that
\[
\bar\pa_0\nabla_0^{1,0}\lb \lb 1+\bar\Delta_{f_0}\rb\eta_{i\bar j}-\tr\lb\phi_i\bar\phi_j\rb\rb=0.
\]
Since $\lfb\omega_t\rfb$ is $T$-invariant, we know that $\rho$, and hence also $\eta_{i\bar j}$, are $T$-invariant. This implies that $\bar\Delta_{f_0}\eta_{i\bar j}=\Delta_{f_0}\eta_{i\bar j}$ and then formula \eqref{ip100} follows from Theorem \ref{e1}.

Finally, by formula \eqref{ip100}, we know that
\[
\int_{M_0}\lb \lb 1+\Delta_{f_0}\rb\eta_{i\bar j}-\tr\lb\phi_i\bar\phi_j\rb\rb\Omega_0=0
\]
and thus
\begin{align*}
\begin{split}
\int_{M_0}\eta_{i\bar j}\ \Omega_0= & \int_{M_0}\lb 1+\Delta_{f_0}\rb\eta_{i\bar j}\ \Omega_0\\
=& \int_{M_0}\tr\lb\phi_i\bar\phi_j \rb\Omega_0\\
= & \left\langle \frac{\pa}{\pa t_i}, \frac{\pa}{\pa t_j}\right\rangle_{WP}(0).
\end{split}
\end{align*}

\end{proof}

\noindent{\it Remark 3.5.} Theorem~\ref{2nd} indicates the close relation between the variation of \kr solitons and the \wpm metric we have defined. By finding higher order expansions of the function $\rho$ defined by \eqref{ip35},  we can derive the curvature formula of the \wpm metric. However, note that there is no apparent bound of the curvature of the \wpm metric even in the case of Fano \ke manifolds.

\medskip
\bigskip
\noindent {\bf Declaration}

 \bigskip
\noindent {\bf Conflict of interest}  \ There is no conflict of interest to disclose.


\begin{thebibliography}{99}

\bibitem{ahlforsa}
L. Ahlfors,
\newblock Some remarks on Teichm\"uller's  space of Riemann surfaces.
\newblock {\em Ann. of Math. (2) \textbf {74} (1961), 171--191.}

\bibitem{ahlforsb}
L. Ahlfors,
 \newblock Curvature properties of Teichm\"uller's space.
 \newblock {\em J. Analyse Math.  {\textbf 9} (1961/1962), 161--176.}

\bibitem{cgb88}
P.~Candelas, P.~S. Green, and T.~H\"{u}bsch,
\newblock Connected {C}alabi-{Y}au compactifications (other worlds are just
  around the corner).
\newblock In {\em Strings'88} ({C}ollege {P}ark, {MD}, 1988), 155--190.
  World Sci. Publ., Teaneck, NJ, 1989.

\bibitem{cao} H.-D. Cao, Existence of gradient K\"{a}hler-Ricci solitons, {\em Elliptic and Parabolic Methods in Geometry} (Minneapolis, MN, 1994), A K Peters, Wellesley, MA, (1996) 1--16.

\bibitem{cao2010} H.-D. Cao,
\newblock Recent progress on Ricci solitons,
\newblock  {\em Recent advances in geometric analysis}, 1--38, Adv. Lect. Math. (ALM), {\bf 11} Int. Press, Somerville, MA, 2010.


\bibitem{csyz2021}
H.-D. Cao, X.~Sun, S.-T. Yau, and Y.~Zhang,
\newblock On the deformation of Fano manifolds.
\newblock {\em Math. Ann.}  (2021). https://doi.org/10.1007/s00208-021-02226-2

\bibitem{chu}T. Chu, ``The Weil-Petersson metric in the moduli space", PhD thesis, Columbia University, 1976.

\bibitem{don97}S. K. Donaldson, Remarks on gauge theory, complex geometry and 4-manifold topology.  \emph{Fields Medallists' lectures}, 384-403, World Sci. Ser. 20th Century Math., 5, World Sci. Publ., River Edge, NJ, 1997.

 \bibitem{do2015}
S. K. Donaldson,
\newblock The Ding functional, Berndtsson convexity and moment maps.
\newblock {\em Geometry, analysis and probability, 57-67, Progr. Math., 310, Birkhäuser/Springer, Cham, 2017}.

\bibitem{fs1990}
A.~Fujiki and G.~Schumacher,
\newblock The moduli space of extremal compact {K}\"{a}hler manifolds and
  generalized {W}eil-{P}etersson metrics.
\newblock {\em Publ. Res. Inst. Math. Sci.}, \textbf {26}(1):101--183, 1990.

\bibitem{fuj}A. Fujiki, Moduli space of polarized algebraic manifolds and Kähler metrics,  \emph{Sugaku Expositions}, {\bf 5} (1992), no. 2, 173--191.

\bibitem{futaki}
A.~Futaki,
\newblock K\"ahler-{E}instein metrics and integral invariants, Vol. \textbf {1314},  {\em Lecture Notes in Mathematics}.
\newblock Springer-Verlag, Berlin, 1988.

\bibitem{he}
W.~He,
\newblock K\"ahler-Ricci soliton and H-functional,
\newblock {\em Asian J. Math.} {\textbf 20} (2016), no. 4, 645--663.

\bibitem{koiso83}N. Koiso, Einstein metrics and complex structures.  \emph{Invent. Math}. {\bf 73} (1983), no. 1, 71-106.

\bibitem{koiso90} Koiso, N., On rotationally symmetric Hamilton's equation for K\"ahler-Einstein metrics.{\em  Recent topics in differential and analytic geometry}, 327--337, Adv. Stud. Pure Math., \textbf{18-I}, Academic Press, Boston, MA, 1990.

\bibitem{macama}
T.~Mabuchi.
\newblock A theorem of {C}alabi-{M}atsushima's type.
\newblock {\em Osaka J. Math.}, 39(1):49--57, 2002.

\bibitem{mas}H. Masur, The extension of the Weil-Petersson metric to the
boundary of Teichmuller space.  \emph{Duke Math. J}. {\bf 43} (1976), 623--635.

\bibitem{nan}
A.~Nannicini,
\newblock Weil-{P}etersson metric in the moduli space of compact polarized
  {K}\"{a}hler-{E}instein manifolds of zero first {C}hern class.
\newblock {\em Manuscripta Math.}, {\bf 54} (1986), no. 4, 405--438.




\bibitem{prsa21}
O. Garc\'ia-Prada and D. Salamon,
\newblock A moment map interpretation of the Ricci form, K\"ahler-{E}instein structures, and Teichmüller spaces.
\newblock {\em Integrability, quantization, and geometry II. Quantum theories and algebraic geometry, 223–255, Proc. Sympos. Pure Math., 103.2, Amer. Math. Soc., Providence, RI, [2021], ©2021.}


\bibitem{ro1974}
H.L. Royden,
\newblock Intrinsic metrics on Teichm\"uller space.
\newblock {Proceedings of the International Congress of Mathematicians (Vancouver, B. C., 1974), Vol. 2, pp. 217–221. Canad. Math. Congress, Montreal, Que., 1975. }

\bibitem{sch}
G. Schumacher,
\newblock The curvature of the Petersson-Weil metric on the moduli space of \ke manifolds.
\newblock {\em Complex analysis and geometry, 339-354, Univ. Ser. Math., Plenum, New York, 1993.}

\bibitem{siu}
Y.-T. Siu,
\newblock Curvature of the Weil-Petersson metric in the moduli space of compact \ke manifolds of negative first Chern class.
\newblock {\em Contributions to several complex variables, 261-298,
Aspects Math., E9, Friedr. Vieweg, Braunschweig, 1986.}


\bibitem{sundeform1}
X.~Sun,
\newblock Deformation of canonical metrics {I}.
\newblock {\em Asian J. Math.}, \textbf {16} (2012), no.1, 141--155.

\bibitem{tod}A. Todorov, The Weil-Petersson geometry of the moduli space of SU($n\ge 3$) (Calabi-Yau) manifolds. I.  \emph{Comm. Math. Phys.} {\bf 126} (1989), no. 2, 325--346.

\bibitem{tianzhu}
G. Tian and X. H. Zhu,
\newblock A new holomorphic invariant and uniqueness of K\"ahler-Ricci solitons.
\newblock {\em Comment. Math. Helv}. \textbf {77} (2002), no. 2, 297--325.

\bibitem{trwil11}
T.~Trenner and P.~M.~H. Wilson,
\newblock Asymptotic curvature of moduli spaces for {C}alabi-{Y}au threefolds.
\newblock {\em J. Geom. Anal.} {\bf 21} (2011), 409--428.

\bibitem{wangzhu} X.J. Wang and X.H. Zhu, {K$\ddot{a}$hler-Ricci solitons on toric manifolds with
positive first Chern class}, {\em Adv. Math.} \textbf{188 }(2004), no.1, 87--103.

\bibitem{weil} A. Weil, Modules des surfaces de Riemann. (French),  \emph{1958 S\'eminaire Bourbaki}; 10e ann\'ee: 1957/1958. Textes des conf\'erences; Expos\'es 152\`a 168; 2e \'ed.corrig\'e, Expos\'e 168 7 pp. Secr\'etariat math\'ematique, Paris.

\bibitem{wilson04}
P.~M.~H. Wilson,
\newblock Sectional curvatures of {K}\"{a}hler moduli.
\newblock {\em Math. Ann.} {\bf 330} (2004), 631--664.

\bibitem{wolpert}S. Wolpert, Chern forms and the Riemann tensor for the moduli space of curves.  \emph{Invent. Math}. {\bf 85} (1986), no. 1, 119--145.

\bibitem{zhang}
Y.~Zhang,
\newblock ``Geometric quantization of classical metrics on the moduli space of canonical metrics", Ph.D. thesis, Lehigh University (2014).
\newblock https://preserve.lehigh.edu/etd/1689

\end{thebibliography}
\end{document}